\documentclass[12pt]{article}
\usepackage{graphicx}
\usepackage{epsfig}
\usepackage{fancyhdr}

\usepackage{amsmath,amsfonts,amssymb,amsthm}
\theoremstyle{plain}
\newtheorem{thm}{Theorem}
\newtheorem*{thm*}{Theorem}

\newtheorem*{lem*}{Lemma}
\newtheorem*{cor}{Corollary}

\theoremstyle{definition}

\newtheorem{Def}{Definition}

\newcommand{\reftit}{\textit}    
\newcommand{\refis}{\textbf}     

\begin{document}

\title{Orthogonality and  probability: mixing times}
\author{Yevgeniy Kovchegov\footnote{
 Department of Mathematics,  Oregon State University, Corvallis, OR  97331-4605, USA
 \texttt{kovchegy@math.oregonstate.edu}}}
\date{ }
\maketitle

\abstract{We produce the first example of bounding total variation distance to stationarity and estimating mixing times via orthogonal polynomials diagonalization of discrete reversible Markov chains, the Karlin-McGregor approach. }

\section{Introduction}

 If $P$ is a reversible Markov chain over a sample space $\Omega$, and $\pi$ is a reversibility function (not necessarily a probability distribution), then $P$ is a self-adjoint operator in $\ell^2(\pi)$, the space generated by the inner product
 $$<f,g>_{\pi}=\sum_{x \in S} f(x)g(x)\pi(x)$$
 induced by $\pi$. If $P$ is tridiagonal operator (i.e. a nearest-neighbor random walk) on  $\Omega=\{0,1,2,\dots\}$, then it must have a simple spectrum, and is diagonalizable via orthogonal polynomials as it was studied in the 50's and 60's by Karlin and McGregor, see \cite{km2a}, \cite{szego}. There the extended eigenfuctions $Q_j(\lambda)$ ($Q_0 \equiv 1$) are orthogonal polynomials with respect to a probability measure $\psi$ and
 $$p_t(i,j)=\pi_j \int_{-1}^1 \lambda^t Q_i(\lambda) Q_j(\lambda) d\psi(\lambda)~~~\forall i,j \in \Omega, $$
 where $\pi_j$ ($\pi_0=1$) is the reversibility measure of $P$.
 
 In this paper we are testing a possibility of calculating mixing rates using Karlin-McGregor diagonalization with orthogonal polynomials.
  In order to measure the rate of convergence to a stationary distribution, the following distance is used. 
 \begin{Def}
If $\mu$ and $\nu$ are two probability distributions over a sample space $\Omega$, then the {\it total variation distance} is 
$$\| \nu - \mu \|_{TV} = {1 \over 2} \sum_{x \in \Omega} |\nu(x)-\mu(x)|=\sup_{A \subset \Omega} |\nu(A)-\mu(A)|$$
Observe that the total variation distance measures the coincidence between the distributions on a scale from zero to one.
\end{Def}
 \noindent
 If $\rho=\sum_{k=0}^{\infty} \pi_k < \infty$, then $\nu={1 \over \rho}\pi$ is the stationary probability distribution. If in addition, the aperiodic nearest neighbor Markov chain originates at site $i$, then
 the total variation distance between the distribution $\mu_t=\mu_0P^t$ and $\nu$ is given by
 \begin{eqnarray*}
 \left\|\nu - \mu_t \right\|_{TV} 
 & = & {1 \over 2} \sum_{j} \pi_j \left|\int_{(-1,1)} \lambda^t Q_i(\lambda) Q_j(\lambda) d\psi(\lambda)\right|,
 \end{eqnarray*} 
 as measure $\psi$ contains a point mass of weight ${1 \over \rho}$ at $1$,  
 see \cite{kov}.
 
 The rates of convergence are quantified via mixing times. In the case of a Markov chain over an infinite state space with a unique stationary distribution, the notion of a mixing time depends on the state of origination of the chain.
\begin{Def}
 Suppose $P$ is a Markov chain with a stationary probability distribution $\nu$ that commences at $X_0=i$. Given an $\epsilon >0$, the mixing time $t_{mix}(\epsilon)$ is defined as
 $$t_{mix}(\epsilon)=\min\left\{t~:~\|\nu-\mu_t\|_{TV} \leq \epsilon \right\}$$
\end{Def}
\vskip 0.2 in
\noindent
 In the case of a nearest-neighbor process  on $\Omega=\{0,1,2,\dots\}$ commencing at $i$,  the corresponding mixing time has the following simple expression in orthogonal polynomials
   $$t_{mix}(\epsilon)=\min\left\{t~:~\sum_{j} \pi_j \left|\int_{(-1,1)} \lambda^t Q_i(\lambda) Q_j(\lambda) d\psi(\lambda)\right| \leq 2 \epsilon \right\}$$ 
   Observe that the above expression is simplified when $i=0$. Here we concentrate on calculating mixing times for simple positive recurrent nearest-neighbor Markov chains over $\Omega$, originating from $i=0$.  Our main result concerns the distance to stationarity for a simple random walk with a drift.
 In the main theorem and its corollary we will explore the following Markov chain 
 $$P=\left(\begin{array}{ccccc}0 & 1 & 0 & 0 & \dots \\q & r & p & 0 & \dots \\0 & q & r & p & \ddots \\0 & 0 & q & r & \ddots \\\vdots & \vdots & \ddots & \ddots & \ddots\end{array}\right) \qquad q>p,~~~r>0$$
 
\begin{thm} \label{thm1} Suppose the above Markov chain begins at the origin, $i=0$. Consider the orthogonal polynomials $Q_n$ for the chain. Then the integral
$\int_{(-1,1)} \lambda^t Q_n(\lambda) d\psi(\lambda)$ can be expressed as\\
${(1+q-p)(q+r)-q \over (1+q-p)(q+r)}\cdot\left(-{q \over q+r}\right)^{t+n}
+\left(\sqrt{q \over p}\right)^n \left({p \over q+r}\right) {1 \over 2\pi i}\oint_{|z|=1}{(\sqrt{pq}(z+z^{-1})+r)^t z^n (z-z^{-1}) \over \left(z-\sqrt{p \over q}{r+(1+q-p) \over 2(q+r)}\right)\left(z-\sqrt{p \over q}{r-(1+q-p) \over 2(q+r)}\right)}dz$
and the total variation distance $\left\|\nu - \mu_t \right\|_{TV} $ is bounded above by\\
$$A\left({q \over q+r}\right)^t +B(r+2\sqrt{pq})^t,$$
where $A={(1+q-p)(q+r)-q \over (1+q-p)(1-2p)}$ and $B={\left({p \over q+r}\right)\left(1+{1 \over \sqrt{pq}-p}\right) \over \left(1-\sqrt{p \over q}{r+(1+q-p) \over 2(q+r)}\right)\left(1+\sqrt{p \over q}{r-(1+q-p) \over 2(q+r)}\right)}$.
Therefore, taking $\varepsilon \downarrow 0$, the mixing time $$t_{mix}(\varepsilon)=O\left({\log(\varepsilon) \over \log m(p,q)}\right),$$
where $m(p,q)=\max\left[(r+2\sqrt{pq}),\left({q \over q+r}\right)\right]$.
\end{thm}
 Observe that in the above complex integral all three finite poles are located inside the unit circle. Thus we only need to consider a pole at infinity.

The proof is provided in section \ref{proof}. The result in Theorem \ref{thm1} is the first instance the Karlin-McGregor orthogonal polynomials approach is used to estimate mixing rates. As it was suggested in \cite{kov} we would like the approach to work for a larger class of reversible Markov chains over an infinite state space with a unique stationary distribution. 
There is an immediate corollary (see section \ref{proof}):
\begin{cor}
 If ${q \over q+r}>r+2\sqrt{pq}$, 
 $$\left\|\nu - \mu_t \right\|_{TV} \geq A\left({q \over q+r}\right)^t - B(r+2\sqrt{pq})^t$$
for $t$ large enough, i.e. we have a lower bound of matching order. 
\end{cor}
 Observe that one can easily adjust these results for any origination site $X_0=i$.
 In the next section we will compare the above Karlin-McGregor approach to some of the classical techniques for estimating the ``distance to stationarity" $\left\|\nu - \mu_t \right\|_{TV}$.

 \section{Comparison to the other techniques}

 For the case of geometrically ergodic Markov chains, there are several techniques that produce an upper bound on the distance to stationarity that were developed specifically for the cases when the sample space is large, but finite. These methods are not directly applicable to chains on general state spaces. The coupling method stands out as the most universal. Here we compare the geometric rate in Theorem \ref{thm1} to the one obtained via a classical coupling argument. Then we explain why other geometric ergodicity methods based on renewal theory will not do better than coupling.  
See \cite{mt} and \cite{lindvall} for detailed overview of geometric convergence and coupling.

 \subsection{Geometric convergence via coupling}
 Consider a coupling process $(X_t,Y_t)$, where $X_0=0$ as in Theorem \ref{thm1}, while $Y_0$ is distributed according to the stationary distribution $\nu={1 \over \rho}\pi$. A classical Markovian coupling construction allows $X_t$ and $Y_t$ evolve independently until the coupling time $\tau_{coupling}=\min\{t:~X_t=Y_t\}$.  It is natural to compare $P(\tau_{coupling}>t)$ to
 $P(\tau>t)$, where $\tau=\min\{t:~Y_t=0\}$ is a hitting time, as the chain is positive recurrent.

Now, simple combinatorics implies, for $k \geq n$,
$$P(\tau=k~|~Y_0=n)=\sum_{i,j:~2i+j=k-n} {k! \over i!(i+n)!j!} p^iq^{i+n}r^j$$
Therefore
$$P(\tau>t) \leq {1 \over p\rho} \sum_{k:~k>t}\left( \sum_{i,j,n:~2i+j=k-n} {k! \over i!(i+n)!j!} p^{i+n}q^ir^j \right),$$
where `$\leq$' appears because $\pi_0=1<{1 \over p}$, but it does not change the asymptotic rate of convergence, i.e. we could write `$\approx$' instead of `$\leq$'.
The right hand side can be rewritten as
$${1 \over p\rho} \sum_{k:~k>t} \sum_{j=0}^k \left(\begin{array}{c}k \\j\end{array}\right)r^j(p+q)^{k-j}\ell(k-j)$$ 
for $\ell(m)=P(Y \geq m/2)$, where $Y$ is a binomial random variable with parameters $\left(m,\widetilde{p}={p \over p+q} \right)$.
Now, by Cram\'er's theorem, $\ell(m) \sim e^{[\log2+{1 \over 2}\log\widetilde{p}+{1 \over 2}\log(1-\widetilde{p})]m}$, and therefore
\begin{equation} \label{tau}
P(\tau>t) \sim {1 \over p\rho} \sum_{k:~k>t}(r+2\sqrt{pq})^k={r+2\sqrt{pq} \over p\rho(\sqrt{q}-\sqrt{p})^2}(r+2\sqrt{pq})^t
\end{equation}
Recall that in the Corollary, if $q$ is sufficiently larger than $r$ and $p$, then $\left({q \over q+r}\right)^t$ dominates 
$(r+2\sqrt{pq})^t$, and the total variation distance  $$\left\|\nu - \mu_t \right\|_{TV}=A\left({q \over q+r}\right)^{t} \pm B(r+2\sqrt{pq})^t,$$
where $A$ and $B$ are given in Theorem \ref{thm1} of this paper.
 Thus we need to explain why, when $q$ is sufficiently large, in the equation (\ref{tau}), we fail to notice the dominating term of $\left({q \over q+r}\right)^t$. In order to understand why, observe that the second largest eigenvalue $\left(- {q \over q+r} \right)$ originates from the difference between $\tau_{coupling}$ and $\tau$. In fact, $Y_t$ can reach state zero without ever sharing a site with $X_t$ (they will cross each other, of course).  Consider the case when $p$ is either zero, or close to zero. There, the problem reduces to essentially coupling a two state Markov chain with transition probabilities $p(0,1)=1$ and $p(1,0)={q \over q+r}$. Thus the coupling time will be expressed via a geometric random variable with the failure probability of ${q \over q+r}$.
 
 Of course, one could make the Markov chain $P$  ``lazier" by increasing $r$ at the expense of $p$ and $q$, while keeping proportion ${q \over p}$ fixed, i.e. we can consider $P_{\varepsilon}={1 \over 1+\varepsilon}(P+\varepsilon I)$. This will minimize the chance of $X_t$ and $Y_t$ missing each other, but this also means increasing $(r+2\sqrt{pq})$, and slowing down the rate of convergence in (\ref{tau}).

 In order to obtain the correct exponents of convergence, we need to {\bf redo the coupling rules} as follows. We now let the movements of $X_t$ and $Y_t$ be synchronized whenever both are not at zero (i.e. $\{X_t,Y_t\} \cap \{0\} =\emptyset$), while letting  $X_t$ and $Y_t$ move independently when one of them is at zero, and the other is not. Then at the hitting time $\tau$,
 either $X_t=Y_t=0$ and the processes are successfully coupled, or $X_t=1$ and $Y_t=0$. In the latter case we are back to the geometric variable with the failure probability of ${q \over q+r}$. That is, the only way for $X_t$ and $Y_t$ to couple would be if one of the two is at state $0$ and  the other is at state $1$. Using the set theory notations, if $\{X_t,Y_t\}=\{0,1\}$, conditioning on $\{X_{t+1},Y_{t+1}\} \not= \{1,2\}$ would give us
$$\{X_{t+1},Y_{t+1}\}=
 \begin{cases}
      \{1\} & \text{ with probability }{r \over q+r}, \\
      \{0,1\} & \text{ with probability }{q \over q+r},
\end{cases}$$
 
 When ${q \over q+r}>r+2\sqrt{pq}$, the above modified coupling captures the order $\left({q \over q+r}\right)^t$. The coefficient $A$ however is much harder to estimate using the coupling approach, while it is immediately provided in Theorem \ref{thm1} and its corollary. Take for example $p={1 \over 11}$,  $r={1 \over 11}$ and $q={9 \over 11}$. There  ${q \over q+r}>r+2\sqrt{pq}$, and according to the Corollary, the lower bound of $A\left({q \over q+r}\right)^{t} - B(r+2\sqrt{pq})^t$ and the upper bound of $A\left({q \over q+r}\right)^{t} + B(r+2\sqrt{pq})^t$ are of the matching order, and the oreder of convergence is tight
  $$\left\|\nu - \mu_t \right\|_{TV}={91 \over 171}\left({9 \over 10}\right)^t \pm {39 \over 28} \left({7 \over 11}\right)^t $$ 

 \subsection{Drift, minorization and geometric ergodicity}
 
  The optimal ``energy function" $V(x)=\left({q \over p}\right)^{x/2}$ converts the geometric drift inequality in Meyn and Tweedie \cite{mt} Chapter 15  into equality
  $$E[V(X_{t+1})~|~X_t=x]=(r+2\sqrt{pq})V(x)+\left(\sqrt{q \over p}- (r+2\sqrt{pq}) \right) 1\!\!1_C(x)$$
  thus confirming the geometric convergence rate of $(r+2\sqrt{pq})^t$ for the tail probability $P(\tau_C>t)$, where $C=\{0\}$ is the obvious choice for the ``small set", and $\tau_C$ is the hitting time. Once again all the challenge is at the origin. 
  In fact there is only a trivial ``minorization condition" when  $C=\{0\}$. The minorization condition reads
  $$p(x,A) \geq \epsilon Q(A)~~~\forall x \in C,~A \subset \Omega,$$
  where, if $C=\{0\}$,  the only choice for the probability measure $Q$ is $Q=\delta_1$, and $\epsilon=1$.
  With $\epsilon=1$ the split of the Markov chain is trivial, and as far as the corresponding coupling goes, the only issue would be (as we mentioned before) to compute  the tail of the hitting time $\min\{t:~(X_t,Y_t) \in C \times C\}$ when $q$ is large. 
  If $C=\{0,1,\dots,k\}$ for some $k >0$, there is no minorization condition. In the latter case, estimating the hitting time $\min\{t:~(X_t,Y_t) \in C \times C\}$ is straightforward,  but without minorization, this will not be enough to estimate the tail for the coupling time. The ``splitting technique" will not work,  rather a coupling approach of the preceding subsection to be pursued.
  
 The case of recurrent reflecting random walk (the M/M/1 queue) had been considered as one of the four benchmark examples in the geometric ergodicity theory (see \cite{bax} and references therein). There, in the absence of the second largest eigenvalue of $\left(-{q \over q+r}\right)$, with $r=0$,  the rate of $(2\sqrt{pq})^t$ was proven to be the optimal (see \cite{lund}).
 The methods in the theory of geometric convergence are in the most part based on the renewal theory (see \cite{mt}, \cite{bax}, \cite{rose} and references therein), and concentrate more on the tail probability for the hitting time $\tau_C$ in the splitting method.
 As for the Markov chain $P$ considered in this paper, in the absence of a useful splitting of it, the approach that works is the coupling. While the coupling provides the right exponents, it does not necessarily produce the tight coefficients.


 \section{The proof of Theorem \ref{thm1}} \label{proof}
 
 \begin{proof}
 Since we require $r>0$ for aperiodicity, we will need to obtain the spectral measure $\psi$ via an argument similar to that of Karlin and McGregor in \cite{km2a}, where the case of $r=0$ was solved. 
The orthogonal polynomials are obtained via solving a simple linear recursion:
 $Q_0=1$, $Q_1=\lambda$, and
$Q_n(\lambda)=c_1(\lambda) \rho^n_1(\lambda)+c_2(\lambda)\rho^n_2(\lambda)$,
 where  $\rho_1(\lambda)={\lambda-r+\sqrt{(\lambda-r)^2-4pq} \over 2p}$ and 
 $\rho_2(\lambda)={\lambda-r-\sqrt{(\lambda-r)^2-4pq} \over 2p}$
  are the roots of the characteristic equation for the recursion,
  and $c_1={\rho_2-\lambda \over \rho_2 -\rho_1}$ and $c_2={\lambda-\rho_1 \over \rho_2 -\rho_1}$.
  
\vskip 0.2 in
\noindent  
 Now $\pi_0=1$, $\pi_n={p^{n-1} \over q^n}$ ($n\geq 1$) and $\rho={q-p+1 \over q-p}$.  Also, we observe that 
 $$
 \begin{cases}
      |\rho_2(\lambda)| >\sqrt{q \over p} & \text{ on }[-1,r-2\sqrt{pq}), \\
      |\rho_2(\lambda)| <\sqrt{q \over p} & \text{ on }(r+2\sqrt{pq}, 1], \\
      |\rho_2(\lambda)|=\sqrt{q \over p} & \text{ on }[r-2\sqrt{pq}, r+2\sqrt{pq}],
\end{cases}
$$
 and $\rho_1\rho_2={q \over p}$.\\
 The above will help us to identify the point mass locations in the measure $\psi$ since each point mass in $\psi$ occurs when 
 $\sum_{k} \pi_k Q_k^2(\lambda) < \infty$. Thus we need to find all $\lambda \in (r+2\sqrt{pq}, 1]$ such that
 $c_1(\lambda)=0$ and all $\lambda \in [-1,r-2\sqrt{pq})$ such that $c_2(\lambda)=0$. There are two roots, $\lambda=1$ and $\lambda=-{q \over q+r}$.  
 \vskip 0.2 in
\noindent
 We already know everything about the point mass at $\lambda=1$: $Q_k(1)=1$ for all $k \geq 0$, and $\rho=\sum_{k=0}^{\infty} \pi_k Q_k^2(1)={1+q-p \over q-p}$ is the reciprocal  of the point mass at $\lambda=1$.
 \vskip 0.2 in
\noindent
 The only other point mass is at $\lambda=-{q \over q+r}$. One can verify that $\rho_1\left(-{q \over q+r}\right)=-{q \over q+r}$ and $Q_k\left(-{q \over q+r}\right)=\left(-{q \over q+r}\right)^k$, and therefore
$$\sum_{k=0}^{\infty} \pi_k Q_k^2\left(-{q \over q+r}\right)=1+{q \over (q+r)^2-pq}={(1+q-p)(q+r) \over (1+q-p)(q+r)-q}$$
is the reciprocal  of the point mass at $\lambda=-{q \over q+r}$.
\vskip 0.2 in
\noindent
 It follows that the rest of the mass of $\psi$ (other than the two point masses) is spread inside $[r-2\sqrt{pq}, r+2\sqrt{pq}]$. In order to find the density of $\psi$ inside $[r-2\sqrt{pq}, r+2\sqrt{pq}]$ we need to find $(e_0,(P-sI)^{-1}e_0)$ for $Im(s) \not=0$, i.e. the upper left element in the resolvent of $P$.  Let  $(a_0(s),a_1(s),\dots)^T=(P-sI)^{-1}e_0$, then
 $$-sa_0+a_1=1,~~~\text{ and }~~~qa_{n-1}+(r-s)a_n+pa_{n+1}=0$$
 Thus $a_n(s)=\alpha_1 \rho_1(s)^n+\alpha_2 \rho_2(s)^n,$
  where $\alpha_1={a_0(\rho_2-s) -1 \over \rho_2(s)-\rho_1(s)}$ and $\alpha_2={1-a_0(\rho_1-s) \over \rho_2(s)-\rho_1(s)}$.
  \vskip 0.2 in
  \noindent
  Since $(a_0,a_1,\dots) \in \ell^2(\mathbb{C}, \pi)$,
  $$|a_n| \sqrt{q^n \over p^n} \rightarrow 0 \qquad \text{ as } ~~n \rightarrow +\infty$$
 Hence when  $|\rho_1(s)| \not= |\rho_2(s)|$, either $\alpha_1=0$ or $\alpha_2=0$, and therefore
 \begin{equation}\label{a0}
 a_0(s)={\chi_{|\rho_1(s)|< \sqrt{q \over p}} \over \rho_1(s)-s}+{\chi_{|\rho_2(s)|<\sqrt{q \over p}} \over \rho_2(s)-s}
 \end{equation}
\vskip 0.2 in
\noindent
 Now, because of the point masses at $1 $ and $-{q \over q+r}$, $a_0(s)=\int_{(-1,1]}{d\psi(z) \over z-s}$ can be expressed as
$$a_0(s)={q-p \over 1+q-p}\left({1 \over 1-s}\right)+{(1+q-p)(q+r)-q \over (1+q-p)(q+r)}\left({1 \over -{q \over q+r}-s}\right)+\int_{(-1,1)}{\varphi(z)dz \over z-s},$$
where $\varphi(z)$ is an atom-less function. Next we will use the following basic property of Cauchy transforms $Cf(s)={1 \over 2\pi i} \int_{\mathbb{R}}{f(z)dz \over z-s}$ that  
can be derived using the Cauchy integral formula, or similarly, an approximation to the identity formula
\footnote{The curve in the integral does  not need to be $\mathbb{R}$ for $C_+-C_-=I$ to hold.}:
\begin{equation} \label{cauchy}
C_+-C_-=I
\end{equation} 
Here $C_+f(z)=\lim_{s \rightarrow z:~ Im(s)>0} Cf(s)$ and $C_-f(z)=\lim_{s \rightarrow z:~ Im(s)<0} Cf(s)$ for all $z \in \mathbb{R}$.
The above equation (\ref{cauchy}) implies
$$\varphi(x)={ 1\over 2\pi i} \left(\lim_{s=x+i\varepsilon~:~\varepsilon \rightarrow 0+} a_0(s) -\lim_{s=x-i\varepsilon~:~\varepsilon \rightarrow 0+} a_0(s) \right)$$
for all $x \in (-1,1)$.
Recalling (\ref{a0}), we express $\varphi$ as
$\varphi(x)={\rho_1(x)-\rho_2(x) \over  2\pi i(\rho_1(x)-x)(\rho_2(x)-x)}$
 for  $x \in (r-2\sqrt{pq}, r+2\sqrt{pq})$, which in turn simplifies to
 $$\varphi(x)=\begin{cases} 
      {\sqrt{(x-r)^2-4pq} \over 2\pi i ((r+q)x+q)(1-x)}  & \text{ if } x \in (r-2\sqrt{pq}, r+2\sqrt{pq}), \\
        0 & \text{ otherwise }
\end{cases}$$
Let $\mathcal{I}=(r-2\sqrt{pq}, r+2\sqrt{pq})$ denote the support interval, and let $1\!\!1_{I}(x)$ be its indicator function.
Here $$\int_{-1}^1 \varphi(x) dx={p \over q+r}$$
and one can check that
$$\psi(x)={q-p \over 1+q-p}\cdot \delta_1(x)+{(1+q-p)(q+r)-q \over (1+q-p)(q+r)}\cdot \delta_{-{q \over q+r}}(x)+{\sqrt{4pq-(x-r)^2} \over 2\pi ((r+q)x+q)(1-x)} \cdot 1\!\!1_{\cal I}(x)$$
integrates to one.
\vskip 0.2 in
\noindent
Observe that the residues of $g(z)={\sqrt{(z-r)^2-4pq} \over ((r+q)z+q)(1-z)}$ are
$$Res(g(z),1)={q-p \over 1+q-p}~~~\text{ and }~~~
Res\left(g(z),-{q \over q+r}\right)={(1+q-p)(q+r)-q \over (1+q-p)(q+r)}$$
in the principle  branch of the $\log$ function.
\vskip 0.2 in
\noindent
Now\\
$\int_{(-1,1)} \lambda^t Q_n(\lambda) d\psi(\lambda)={(1+q-p)(q+r)-q \over (1+q-p)(q+r)}\cdot\left(-{q \over q+r}\right)^{t+n}
+\int_{r-2\sqrt{pq}}^{r+2\sqrt{pq}} \lambda^t (c_1\rho_1^n+c_2\rho_2^n){\rho_1 -\rho_2 \over 2\pi i(\rho_1- \lambda)(\rho_2-\lambda)} d\lambda$
\vskip 0.16 in
\noindent
and therefore, since $c_1={\rho_2-\lambda \over \rho_2 -\rho_1}$ and $c_2={\lambda-\rho_1 \over \rho_2 -\rho_1}$,
\vskip 0.16 in
$\int_{(-1,1)} \lambda^t Q_n(\lambda) d\psi(\lambda)={(1+q-p)(q+r)-q \over (1+q-p)(q+r)}\cdot\left(-{q \over q+r}\right)^{t+n}
+{1 \over 2\pi i}\int_{r-2\sqrt{pq}}^{r+2\sqrt{pq}} \lambda^t \left({\rho_2^n \over \rho_2-\lambda}-{\rho_1^n \over \rho_1-\lambda}\right)d\lambda,$
\vskip 0.16 in
\noindent
where, if we
let $\rho_1=\sqrt{q \over p} z$ for $z$ in the lower semicircle and $\rho_2=\sqrt{q \over p} z$ for $z$ in the upper semicircle, then
\vskip 0.2 in
${1 \over 2\pi i}\int_{r-2\sqrt{pq}}^{r+2\sqrt{pq}} \lambda^t \left({\rho_2^n \over \rho_2-\lambda}-{\rho_1^n \over \rho_1-\lambda}\right)d\lambda
=\left(\sqrt{q \over p}\right)^n {1 \over 2\pi i}\oint_{|z|=1}{(\sqrt{pq}(z+z^{-1})+r)^t z^n \sqrt{pq}(1-z^{-2})dz\over \sqrt{q \over p}z-(\sqrt{pq}(z+z^{-1})+r)}$
$$=\left(\sqrt{q \over p}\right)^n \left({p \over q+r}\right) {1 \over 2\pi i}\oint_{|z|=1}{(\sqrt{pq}(z+z^{-1})+r)^t z^n (z-z^{-1}) \over \left(z-\sqrt{p \over q}{r+(1+q-p) \over 2(q+r)}\right)\left(z-\sqrt{p \over q}{r-(1+q-p) \over 2(q+r)}\right)}dz$$
Here the absolute value of the function in the last integral is bounded by $M(r+2\sqrt{pq})^t$ with
$M={2 \over \left(1-\sqrt{p \over q}{r+(1+q-p) \over 2(q+r)}\right)\left(1+\sqrt{p \over q}{r-(1+q-p) \over 2(q+r)}\right)}$.
Therefore, plugging in the values of $\pi_n$, we show that the distance to stationarity $\left\|\nu - \mu_t \right\|_{TV}={1 \over 2}\sum_{n=0}^{\infty} \pi_n \left|\int_{(-1,1)} \lambda^t Q_n(\lambda) d\psi(\lambda)\right|$ is bounded above by
$$A \left({q \over q+r}\right)^t+B(r+2\sqrt{pq})^t$$
where  
$$A={(1+q-p)(q+r)-q \over 2(1+q-p)(q+r)} \sum_{n=0}^{\infty}\pi_n \left({q \over q+r}\right)^n={(1+q-p)(q+r)-q \over (1+q-p)(1-2p)}$$
and
$$B={M \over 2} \left({p \over q+r}\right)\left(1+{1 \over \sqrt{pq}-p}\right)={\left({p \over q+r}\right)\left(1+{1 \over \sqrt{pq}-p}\right) \over \left(1-\sqrt{p \over q}{r+(1+q-p) \over 2(q+r)}\right)\left(1+\sqrt{p \over q}{r-(1+q-p) \over 2(q+r)}\right)}$$
The above upper bound can be improved if one obtains a better estimate of the trigonometric integrals involved in the sum. 
\vskip 0.2 in
\noindent
We conclude that
$t_{mix}(\varepsilon)=O\left({\log(\varepsilon) \over \log m(p,q)}\right)$
as $\varepsilon \downarrow 0$.
\end{proof}

\section*{Acknowledgment}
The author would like to acknowledge useful comments about the idea of using orthogonal polynomials for computing mixing times he received from R.Burton,  A.Dembo, P.Diaconis, M.Ossiander, E.Thomann, E.Waymire and J.Zu\~{n}iga.

\bibliographystyle{amsplain}

\end{document}